\documentclass[preprint,12pt]{elsarticle}

\usepackage{amssymb}
\usepackage{pifont}
\usepackage{mathrsfs}
\usepackage{graphicx,ifpdf}
\usepackage{amsmath}
\usepackage[colorlinks,anchorcolor=blue,linkcolor=blue,citecolor=blue]{hyperref}
\usepackage{amsfonts,amsthm,mathrsfs,bbding}
\usepackage{graphics,multicol}
\usepackage{color,fourier}
\usepackage{float}
\usepackage{enumerate}
\usepackage{booktabs}
\usepackage{epstopdf}
\usepackage{indentfirst}
\usepackage{empheq}
\usepackage{dsfont}
\usepackage{latexsym,bm}
\allowdisplaybreaks[4]
\evensidemargin=\oddsidemargin\topmargin=-1.5cm

\newtheorem{conj}{Conjecture}
\newtheorem{lem}{Lemma}[section]
\newtheorem{prob}{Problem}[section]
\newtheorem{thm}{Theorem}[section]
\newtheorem{cor}{Corollary}[section]
\newtheorem{pro}{Proposition}

\newtheorem{fact}[conj]{Fact}

\newtheorem{claim}{Claim}
\newtheorem{definition}{Definition}

\addtocounter{section}{0}

\usepackage{tikz}

\journal{Discrete Mathematics}

\begin{document}

\begin{frontmatter}



\title{On the multiplicity of the Laplacian eigenvalue $1$ of a graph}


\author[label1,label2]{Fei Wen \corref{cor1}}
\ead{wenfei@lzjtu.edu.cn}
\cortext[cor1]{Corresponding author}

\author[label1]{Kai Wang}
\ead{wangkai202402@126.com}

\author[label1,label2]{Tao Wang}
\ead{taowang815@163.com}

\affiliation[label1]{organization={Institute of Applied Mathematics},
	addressline={Lanzhou Jiaotong
		University},
	city={Lanzhou},
	postcode={730070},
	state={Gansu},
	country={PR China}}

\affiliation[label2]{organization={Gansu Center for Fundamental Research in Complex Systems Analysis and Control},
	addressline={Lanzhou Jiaotong
		University},
	city={Lanzhou},
	postcode={730070},
	state={Gansu},
	country={PR China}}

%% Abstract
\begin{abstract}
A connected graph $G$ with $n\ge 2$ vertices and $\varepsilon$ edges is called a \emph{$k$-cyclic graph} if $k=\varepsilon-n+1$. The Laplacian matrix is defined as $L(G)=D(G)-A(G)$, where $D(G)$ is the diagonal matrix of vertex degrees and $A(G)$ is the adjacency matrix of $G$. Let $m_{G}(1)$ be the multiplicity of the Laplacian eigenvalue $1$ of $G$. In this paper, we show that for fixed integers $n\ge \big\lceil\frac{3+\sqrt{8k+1}} {2}\big\rceil$ and $k\geq0$, if $G$ is a $k$-cyclic graph of order $n$, then
\begin{equation*}
	\displaystyle m_{G}(1)\leq  n-\big\lceil\frac{3+\sqrt{8k+1}} {2}\big\rceil
\end{equation*}
and moreover, we completely characterize all extremal graphs attained the upper bound.
\end{abstract}

%% Keywords

\end{frontmatter}


\section{Introduction}
All graphs considered in this paper are assumed to be simple and connected. Let $G=(V(G),E(G))$ be a graph of order
$n$ and size $\varepsilon$, where $V(G)$ and $E(G)$ denote the vertex set and edge set of $G$, respectively. The complement of a graph $G$, denoted by $\overline{G}$, is the graph with vertex set $V(\overline{G})=V(G)$, in which two distinct vertices are adjacent if and only if they are nonadjacent in $G$. We denote by $P_{n}$, $C_{n}(n\ge3)$, $S_{n}$, $K_{n}$  the \emph{path}, the \emph{cycle}, the \emph{star}, \emph{complete graph} on $n$ vertices, respectively and by $\omega(G')$ the number of connected components of a disconnected graph $G'$. Let $G-v$ denote the graph obtained from $G$ by deleting the vertex $v$ and all edges incident with it, and $G-tK_{1}$ denote the graph obtained from $G$ by deleting $t$ isolated vertices. We write as $G+k\cdot e$ the graph obtained from $G$ by adding $k$ edges between its $k$ pairs of non-adjacent vertices. Let $G_{1}=(V(G_{1}),E(G_{1}))$ and $G_{2}=(V(G_{2}),E(G_{2}))$ be two graphs. The \emph{union} $G_{1}\cup G_{2}$ is the graph with vertex set $V(G_{1})\cup V(G_{2})$ and edge set $E(G_{1})\cup E(G_{2})$; when $G_{1}$ and $G_{2}$ are vertex-disjoint, their join, denoted by $G_{1}\vee G_{2}$, is obtained from their disjoint union by adding all possible edges between $V(G_{1})$ and $V(G_{2})$; that is,
$V(G_1\vee G_2)=V(G_1)\cup V(G_2),$
$E(G_1\vee G_2)=E(G_1)\cup E(G_2)\cup\{\,uv: u\in V(G_1),\ v\in V(G_2)\,\}.$ A vertex subset of a graph $G$ is called a \emph{clique} (resp. an \emph{independent set}) if every pair of distinct vertices in the subset is adjacent (resp. nonadjacent) in $G$. The \emph{clique number} (resp. \emph{independence number}) is defined as the maximum number of vertices in a clique (resp. an independent set) of $G$, and is denoted by $\beta(G)$ (resp.\ $\alpha(G)$). A \emph{pendant vertex} of a graph $G$ is a vertex of degree $1$, and a quasi-pendant vertex of $G$ is a vertex that is adjacent to at least one pendant vertex. Let $p(G)$ and $q(G)$ (or $p$ and $q$ if there is no confusion) denote the number of pendant vertices and quasi-pendant vertices in $G$, respectively. A connected graph $G$ of order $n$ and size $\varepsilon$ is called a \emph{$k$-cyclic graph} if $k=\varepsilon-n+1$. Especially, if $k=0$ then $G$ is a \emph{tree} of order $n$, and $k=1$ then $G$ is a \emph{unicyclic graph}, and $k=2$ then $G$ is a \emph{bicyclic graph} of order $n$, and so on.

Let $A(G)$ be the adjacency matrix of $G$, and $D(G)$ the degree diagonal matrix of $G$. Then the Laplacian matrix of $G$ is defined as $L(G)=D(G)-A(G)$. Note that $L(G)$ is positive semidefinite with all row sums being equal to zero. Thus, the Laplacian eigenvalues of $G$ can be arranged in
non-increasing order as
\begin{equation*}
	\mu_{1}(G) \geq \mu_{2}(G) \geq \cdots \geq
	\mu_{n-1}(G) \geq \mu_{n}(G)=0.
\end{equation*}
For $1\leq \ell\leq n$, let $s_{\ell}(G)=\sum\limits_{i=1}^{\ell}\mu_{i}(G)$
denote the sum of the $\ell$ largest Laplacian
eigenvalues of $G$. As usual, let $m_{G}(\mu)$ be the multiplicity of a Laplacian eigenvalue $\mu$ of $G$. Especially, if $\mu=1$ then $m_{G}(1)$ denotes the multiplicity of the Laplacian eigenvalue $1$ of $G$. The Laplacian eigenvalues of $G$, together with their multiplicities, constitute the Laplacian spectrum of $G$, denoted by $\operatorname{Spec}_{L}(G)$. In particular, if $G$ has $s$ distinct Laplacian eigenvalues $\mu_{1}>\mu_{2}>\cdots>\mu_{s}$, then $\operatorname{Spec}_{L}(G)=\{\mu_{1}^{m_{1}},\mu_{2}^{m_{2}},\ldots,\mu_{s}^{m_{s}}\}$, where $m_{i}$ is the multiplicity of $\mu_{i}$ in $G$, that is, $m_{G}(\mu_{i})=m_{i}$, and $\sum_{i=1}^{s}m_{i}=n$.

\begin{figure}[H]\hspace{-2.2cm}
	\centering
	\includegraphics[width=6.7in]{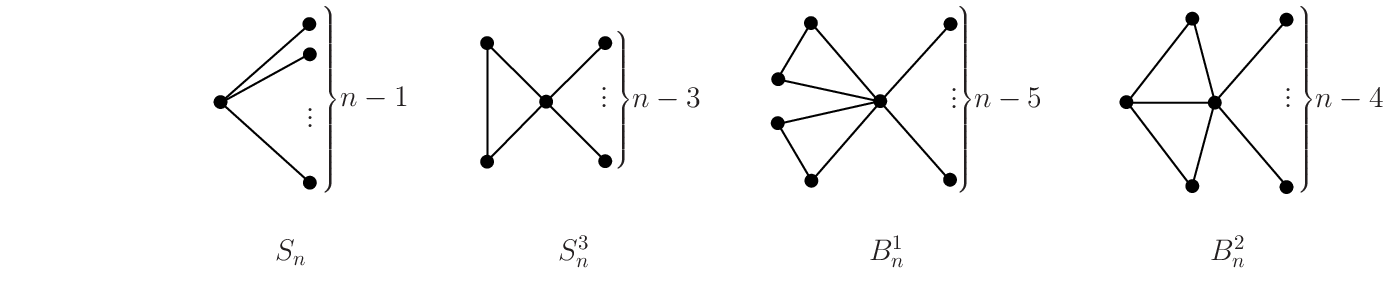}
	\caption{\footnotesize  Graphs $S_{n}$, $S_{n}^{3}$, $B_{n}^{1}$ and $B_{n}^{2}$. }\label{f-1}
\end{figure}

The multiplicity of the Laplacian eigenvalue $1$ has been studied extensively in recent years, owing to its intimate links with various combinatorial invariants and structural properties of graphs, see \cite{Guo,Wen, Wang, Han,Tian1, Tian2, Zhou} and references therein. In 1990, Grone et al. \cite{Grone} established that for a tree $T$ on $n$ vertices, every integral Laplacian eigenvalue $\mu>1$ is simple, i.e., $m_{T}(\mu)=1$, and divides $n$. Subsequently, Guo et al.~\cite{Guo} in 2010 determined the distribution of $m_{T}(1)$ for trees $T$ of order $n\ge 4$, and obtained the following result.
\begin{thm}[\cite{Guo}, Theorem 2.3]\label{thm-k-1}
	Let $T$ be a tree on $n\geq4$ vertices. Then $m_{T}(1)\leq n-2$, where the equality holds if and only if $T\cong S_{n}$, see Fig $\ref{f-1}$.
\end{thm}

In 2022, Wen and Huang \cite{Wen} according to the number of pendant vertices derived a sharp upper bound
\begin{equation}\label{equa-s-1}
	m_{G}(\mu)\leq p(G)+2\varepsilon-2n+1
\end{equation}
for the multiplicity of a Laplacian eigenvalue $\mu$ of a graph $G$, and then they introduced two graph operations, termed \emph{grafting} and \emph{shifting}, and employed them to investigate the distribution of the Laplacian eigenvalue $1$ for unicyclic graphs on $n\ge 7$ vertices; moreover, they characterized all unicyclic graphs attaining the first largest multiplicity of the Laplacian eigenvalue $1$ below.

\begin{thm}[\cite{Wen}, Theorem 4.3]\label{thm-k-2}
	Let $U$ be a unicyclic graph on $n\geq4$ vertices. Then $m_{U}(1)\leq n-3$, where the equality holds if and only if $U\cong S_{n}^{3}$, see Fig $\ref{f-1}$.
\end{thm}
It is worth mentioning that Han et al. in \cite{Han} completely characterized the case where equality holds in (\ref{equa-s-1}) for $\mu=1$.
In recent times, Tian et al. \cite{Tian1, Tian2} showed that
\begin{equation*}
	m_{T}(1)=p(T)-q(T)+m_{\widehat{T}}(1),
\end{equation*}
where $\widehat{T}$ denotes the \emph{reduced tree} of $T$, with $p(\widehat{T})=q(\widehat{T})$. They further obtained upper bounds on the multiplicity of the Laplacian eigenvalue $1$ for reduced trees and unicyclic graphs of order $n\ge10$ without pendant paths $P_{3}$. Additionally, Zhou et al. \cite{Zhou} studied the multiplicity of the Laplacian eigenvalue $1$ of connected \emph{reduced graphs}, and proved that $m_{G}(1)\leq k+\frac{n-2}{4}$ for a connected $k-$cyclic reduced graph $G$ of order $n\geq 6$, where $k$ is the cyclomatic number of $G$, and characterized the extremal graphs achieving equality.

Most recently, Wang et al. \cite{Wang} investigated the distribution of Laplacian eigenvalue $1$ of bicyclic graphs, and characterized  bicyclic graphs with the first largest multiplicity of Laplacian eigenvalue $1$ in the following.

\begin{thm}[\cite{Wang}, Theorem 3.1]\label{thm-k-3}
	Let $B$ be a bicyclic graph on $n\geq5$ vertices. Then $m_{B}(1)\leq n-4$, where the equality holds if and only if $B\cong B_{n}^{1} or B_{n}^{2}$, see Fig $\ref{f-1}$.
\end{thm}

Based on the known results of Theorems \ref{thm-k-1}, \ref{thm-k-2} and \ref{thm-k-3}, the authors \cite{Wang} posed a problem as follows.
\begin{prob}[\cite{Wang}, Problem 1]\label{prob-1}
	Consider the distribution of Laplacian eigenvalue $1$ of connected graphs, with a fixed cyclomatic number, and characterize the extremal graphs attaining the maximum multiplicity of Laplacian eigenvalue $1$.
\end{prob}

The research motivation of this paper arises from two distinct perspectives, the first is the established results outlined above, and the second is rooted in Problem \ref{prob-1}, which asks \emph{which graphs achieve the maximum multiplicity of the Laplacian eigenvalue $1$}. Theorems \ref{thm-k-1}, \ref{thm-k-2} and \ref{thm-k-3} indicate  that a $k$-cyclic graph (i.e., $k=0,1,2$) of order $n\ge 2$, attained the Largest multiplicity of the Laplacian eigenvalue $1$, consists of a $k$-cyclic subgraph with as few vertices as possible together with a number of pendant vertices. Accordingly, we first investigate the minimum order of a $k$-cyclic graph on $n\ge 2$ vertices, leading to the following proposition.
\begin{pro}\label{prop-k-1}
	For a fixed integer $k\geq0$, let $G$ be a $k-$cyclic graph of order $n\geq2$. Then $n\geq \big\lceil\frac{3+\sqrt{8k+1}} {2}\big\rceil$.
\end{pro}
For convenience, we write $n_{k}=\min\big\{n|n\ge \frac{3+\sqrt{8k+1}} {2}\big\}=\big\lceil\frac{3+\sqrt{8k+1}} {2}\big\rceil$. Especially, if $n_{k}=\frac{3+\sqrt{8k+1}}{2}$, we denote by $\widehat{n}_{k}=n_{k}$ for the sake of distinction. A $k$-cyclic graph is called to be a \emph{minimal $k$-cyclic graph} if it has exactly $n_{k}$ vertices. Proposition \ref{prop-k-1} implies that a $k$-cyclic graph is of order at least $n_{k}$. Let $G$ be a $k$-cyclic graph on $n$ vertices. Then it infers from Theorems \ref{thm-k-1}, \ref{thm-k-2} and \ref{thm-k-3} that $m_{G}(1)\leq n-n_{k}$ for $n\ge n_{k}$. Subsequently, we confirm the assertion, and obtain the following result.

\begin{thm}\label{thm-k-5}
	For fixed integers $n\ge n_{k}$ and $k\geq0$, let $G$ be a $k$-cyclic graph of order $n$. Then
		\begin{equation*}\label{equa-w-2}
			\displaystyle m_{G}(1)\leq  n-\big\lceil\frac{3+\sqrt{8k+1}} {2}\big\rceil.
		\end{equation*}
		Moreover, if $n=n_{k}$, then $m_{G}(1)=0$ if and only if $G$ is a minimal $k$-cyclic graph; and otherwise,
		\begin{itemize}
			\item [(i)] if $\displaystyle k={\widehat{n}_{k}-1\choose 2}$, then
			\begin{equation*}\label{equa-w-3}
				\displaystyle m_{G}(1)\leq  n-\frac{3+\sqrt{8k+1}}{2}
			\end{equation*}
			with equality if and only if $G\cong K_{1} \vee(K_{n_{k}-1}\cup(n-n_{k})K_{1})$;
			\item [(ii)] if $\displaystyle{\widehat{n}_{k}-1\choose 2}<k<{\widehat{n}_{k}\choose 2}$, then
			\begin{equation*}\label{equa-w-5}
				\displaystyle m_{G}(1)\leq  n-\big[\frac{3+\sqrt{8k+1}}{2}\big]-1
			\end{equation*}
			with equality if and only if
			\begin{equation*}
				G\cong K_{1} \vee(H\cup(n-n_{k})K_{1})
			\end{equation*}
			where $H=K_{n_{k}-1}-\varepsilon(k)\cdot e$, and $\varepsilon(k)={n_{k}-1\choose 2}-k$. Meanwhile, $K_{1} \vee H$ is a minimal $k$-cyclic graph, and $H$ is of order $n_{k}-1$ and size $k$. In particular, if $\varepsilon(k)=n_{k}-3$, then $G$ is also isomorphic to $K_{1} \vee\big(K_{n_{k}-1}\cup K_{2}\cup(n-n_{k}-2)K_{1}\big)$.
	\end{itemize}
\end{thm}

\section{Proof of Proposition \ref{prop-k-1}}
In this section, we will prove the least order of a $k$-cyclic graph for fixed $k\ge 0$.
\begin{proof}[\textbf{Proof of Proposition \ref{prop-k-1}}]
	Let $G$ be a $k$-cyclic graph of order $n\ge 2$ and size $\varepsilon$. It is well-known that for any connected graph of order $n\ge 2$, its number of edges is at most that of the complete graph with the same order. So we have
	
	\begin{equation}\label{equa-1}
		\varepsilon\leq \frac{n(n-1)}{2}, \ \mbox{and} \ k=\varepsilon-n+1.
	\end{equation}
	Hence, it follows from (\ref{equa-1}) that
	\begin{equation}\label{equa-w-1}
		n^{2}-3n+2-2k\geq 0.
	\end{equation}
	The discriminant of (\ref{equa-w-1}) is given by
	\begin{align*}
		\Delta
		&=(-3)^2-4(2-2k)=1+8k >0 \ \  \ (\mbox{since} \ \ k\geq0).
	\end{align*}
	Thus, the corresponding equation of (\ref{equa-w-1}) has two distinct real roots below,
	\begin{equation*}
		n_{1}=\frac{3-\sqrt{8k+1}}{2}
		\quad\text{and}\quad
		n_{2}=\frac{3+\sqrt{8k+1}}{2}.
	\end{equation*}
	Since the leading coefficient of (\ref{equa-w-1}) is positive, we have that
	\begin{equation*}
		n\leq \frac{3-\sqrt{8k+1}}{2}
		\quad\text{or}\quad
		n\geq \frac{3+\sqrt{8k+1}}{2}.
	\end{equation*}
	Furthermore, if $n\le \frac{3-\sqrt{8k+1}}{2}$ then $n\leq1$ due to $k\geq0$, it is a contradiction.
	Thus, we have $n\geq \big\lceil\frac{3+\sqrt{8k+1}} {2}\big\rceil$, as required.
\end{proof}
\section{Proof of Theorem \ref{thm-k-5}}
In what follows, we present some preliminaries that will be used in the subsequent proof.

\begin{definition}[\cite{Chvátal}, Theorem 1]\label{def-w-1}
	For every graph $G$, the following three conditions are equivalent:
	\begin{itemize}
		\item [$(1)$] $G$ is threshold;
		\item [$(2)$] $G$ has no induced subgraph isomorphic to $2K_{2}$, $P_{4}$ or $C_{4}$;
		\item [$(3)$] there is an ordering $v_{1}, v_{2},\ldots,v_{n}$ of the vertices of $G$ and a partition of $v_{2}, v_{3},\ldots,v_{n}$ into disjoint subsets $P$ and $Q$ such that
		\begin{itemize}
			\item [$\bullet$] every $v_{j}\in P$ is adjacent to all the vertices $v_{i}$ with $i<j$,
			\item [$\bullet$] every $v_{j}\in Q$ is adiacent to none of the vertices $v_{i}$ with $i<j$.
		\end{itemize}
	\end{itemize}
\end{definition}
Clearly, it follows from Definition \ref{def-w-1} that the following fact holds.
\begin{fact}\label{fact-w-1}
	The complement of a threshold graph is also a threshold graph.
\end{fact}

%\begin{lem}[\cite{Wen}, Theorem 3.1]\label{lem-w-1}
%	Let $G$ be a $k$-cycle graph with $p(G)$ pendant vertices. If $G$ is not a cycle, then $m_{G}(\mu)\leq 2k+p(G)-1$, and the equality holds if $T$ is a star $S_{n}$ and $\mu=1$.
%\end{lem}
%\begin{rem}\label{rem-1}
%	The case where the equality in Lemma \ref{lem-w-1} holds with $\mu=1$ has been completely characterized in \cite{Han}.
%\end{rem}

\begin{lem}[\cite{Merris}]\label{lem-2-4}
	Let $G$ be a graph with Laplacian spectrum $\operatorname{Spec}_{L}(G)=\{0=\mu_{n},\mu_{n-1},$ $\ldots,\mu_{2},\mu_{1}\}$. Then the Laplacian spectrum of $\overline{G}$ is presented as $\operatorname{Spec}_{L}(\overline{G})= \{0,n-\mu_{1},n-\mu_{2},\ldots,n-\mu_{n-2},n-\mu_{n-1}\}$, where $\overline{G}$ is the complement of the graph $G$.
\end{lem}

Originally, Brouwer and Haemers \cite{Brouwer} proposed the following conjecture: For a graph $G$ with $1\le \ell \le n-1$, $s_{\ell}(G)\leq \varepsilon+\binom{\ell+1}{2}$. After that, Li and Guo \cite{Li} propose the full Brouwer's Laplacian
spectrum conjecture (see Lemma \ref{lem-2-5}), and proved the conjecture holds for $\ell=2$. Subsequently, Chen et al. in \cite{Chen} presented a concise version of the full Brouwer's conjecture, and then show that the full Brouwer's conjecture holds for complete split graphs as well as two families of their spanning subgraphs, and also holds for $k$-cyclic graphs with $k\in \{0,1,2\}$. More recently, Kothari and Tudose\cite{Kothari} proved the Brouwer's conjecture. In \cite{Cai}, Cai et al. proved that for a simple graph of order $n$ with $\varepsilon$ edges and $1\leq \ell \leq n-1$, $s_{\ell}(G)= \varepsilon+\binom{\ell+1}{2}$ if and only if $G$ is a threshold graph with clique number $\ell+1$, which confirms the full Brouwer conjecture proposed in \cite{Li}.

\begin{lem}[\cite{Li, Cai}, The Full Brouwer's Inequality]\label{lem-2-5}
	For any graph $G$ on $n$ vertices with $\varepsilon$ edges and for each $\ell\in \{1,2,\ldots,n-1\}$, $$\displaystyle s_{\ell}(G)\leq \varepsilon+\binom{\ell+1}{2},$$
	with equality if and only if $G\in \mathcal{G}_{\ell,r,s}(r\geq1,s\geq0)$, where $\mathcal{G}_{\ell,r,s}$ is a family of graphs of order $n=\ell+r+s$ consisting of a clique $K_{\ell}$ and two independent sets $\overline{K_{r}}$ and $\overline{K_{s}}$, such that each vertex in $K_{\ell}$ is adjacent to all vertices in $\overline{K_{r}}$, and for each vertex $v_{i}$ in $\overline{K_{s}}$, $N(v_{i})\subsetneq V(K_{\ell})(i=1,2,\ldots,s)$, and $N(v_{i+1})\subseteq N(v_{i})(i=1,2,\ldots,s-1)$, that is, $G$ is a threshold graph having $n=\ell+r+s$ vertices and clique number $\beta(G)=\ell+1$.
\end{lem}

\begin{lem}[\cite{Kirkland},  Corollary 2.1]\label{lem-2-6}
	Let $G$ be a connected graph on $n \geq 3$ vertices. If $G$ has a cutpoint $v$, then $\mu_{n-1}(G) \leq1$, where equality holds if and only if $v$ is adjacent to every vertex of $G$.
\end{lem}

\begin{lem}[\cite{Faria}]\label{lem-2-7}
	Let $G$ be a graph with $p(G)$ pendant vertices and $q(G)$ quasi-pendant vertices. Then $m_{G}(1)\geq p(G)-q(G)$.
\end{lem}

\begin{lem}[\cite{Grone}, Theorem 2.3]\label{lem-2-8}
	Let $\mu$ be a Laplacian eigenvalue of $T$ for some tree $T$ on $n\geq2$ vertices. Then $m_{T}(\mu)\leq p(T)-1$.
\end{lem}

A cut vertex of a graph $G$ is called a \emph{universal cut vertex} if it is adjacent to all vertices except itself.  Let $M(G)$ be a graph matrix of $G$. Then
\begin{equation*}
	\ker(M)=\{\textbf{z}|M(G)\textbf{z}=0\}
\end{equation*}
is said the \emph{kernel} of $G$.  Meanwhile, let $I_{n}$ denote an \emph{identify matrix} of order $n$, and $\textbf{1}_{n}$ the \emph{all-one vector} of order $n$.
\begin{lem}\label{lem-w-9}
	Let $G$ be a graph on $n\ge 2$ vertices, and $u$ be a universal cut vertex of $G$. Then $m_{G}(1)=\omega(G-u)-1$.
\end{lem}
\begin{proof}
	Let $\zeta_{G}(1)$ be the eigensubspace corresponding to the Laplacian eigenvalue $1$ of $G$. We denote by $G-u=G_{1}\cup G_{2}\cup \ldots \cup G_{s}$, where $G_{i}$ is the component of $G-u$ with $n_{i}$ vertices and $i=1,2,\ldots,s$. Let $V(G)=\{u\}\cup V(G_{1})\cup \ldots \cup V(G_{s})$, and  $\textbf{x}=(x_{1}|\textbf{z}^{T}_{1}|\ldots|\textbf{z}^{T}_{s})^{T}$ be the eigenvector corresponding to the Lapalcian eigenvalue $1$ of $L(G)$, and $x_{1}$ be the entry corresponding to $u$ in $\textbf{x}$, and $\textbf{z}_{i}$ be the sub-block corresponding to the component $G_{i}$ in $\textbf{x}$, where $i=1,2,\ldots,s$. Then $\textbf{x}\bot \textbf{1}_{n}$, it thus follows that  $x_{1}+\sum\limits_{\substack{i=1, v\in V(G_{i})}}^{s} \textbf{z}_{i}(v)=0$.
	Furthermore, from $L(G)\cdot\textbf{x}=1\cdot \textbf{x}$ one can obtain that
	\begin{equation}\label{equa-k-1}
		\begin{bmatrix}
			n-1 & -\textbf{1}_{n_{1}}^{T} & -\textbf{1}_{n_{2}}^{T} & \dots & -\textbf{1}_{n_{s}}^{T} \\ -\textbf{1}_{n_{1}} &
			L(G_{1})+I_{n_{1}}& \textbf{0}_{n_{1}\times n_{2}} & \dots & \textbf{0}_{n_{1}\times n_{s}} \\ -\textbf{1}_{n_{2}} & \textbf{0}_{n_{2}\times n_{1}} & L(G_{2})+I_{n_{2}} & \dots & \textbf{0}_{n_{2}\times n_{s}} \\
			\vdots & \vdots & \ldots & \ddots & \vdots \\
			-\textbf{1}_{n_{s}} & \textbf{0}_{n_{s}\times n_{1}} & \dots  & & L(G_{s})+I_{n_{s}}
		\end{bmatrix}
		\begin{bmatrix}
			x_{1} \\ \textbf{z}_{1} \\ \vdots \\ \textbf{z}_{s}
		\end{bmatrix}
		=
		\begin{bmatrix}
			x_{1} \\ \textbf{z}_{1} \\ \vdots \\ \textbf{z}_{s}
		\end{bmatrix}
	\end{equation}
	From the first row of (\ref{equa-k-1}) we have
	\begin{equation*}
		(n-1)x_{1}-\sum\limits_{\substack{i=1, v\in G_{i}}}^{s} \textbf{z}_{i}(v)=x_{1}
	\end{equation*}
	Note that $x_{1}+\sum\limits_{\substack{i=1,v\in V(G_{i})}}^{s} \textbf{z}_{i}(v)=0$. Thus, we have $x_{1}=0$, and so, $\sum\limits_{\substack{i=1,v\in V(G_{i})}}^{s} \textbf{z}_{i}(v)=0$ and  $\textbf{x}=(0|\textbf{z}_{1}|\ldots|\textbf{z}_{s})^{T}$.
	\\
	From the $i+1$-th row of (\ref{equa-k-1}), one can also obtain
	\begin{equation*}
		L(G_{i})\textbf{z}_{i}+\textbf{z}_{i}-\textbf{0}_{{n_{i}\times 1}}=\textbf{z}_{i} \ \ \ \  \text{for} \ i=1,2,\ldots,s,
	\end{equation*}
	and for each $i$, it follows that
	\begin{equation*}
		L(G_{i})\textbf{z}_{i}=0,
	\end{equation*}
	and so, $\textbf{z}_{i}\in\ker(L(G_{i}))$. Thus, there exist arbitrary constants $l_{i}\in \mathbb{R}$ yielding  $\textbf{z}_{i}=l_{i}\textbf{1}_{n_{i}}$, where $i=1,2,\ldots,s$. In this case,
	one can get
	\begin{equation*}
		\textbf{x}=\begin{bmatrix}
			0 \\ l_{1}\textbf{1}_{n_{1}} \\ \vdots \\ l_{s}\textbf{1}_{n_{s}}
		\end{bmatrix}.
	\end{equation*}
	Since $\textbf{x}\bot \textbf{1}_{n}$ we have a linear equation
		\begin{equation}\label{equ-w-4}
			l_{1}n_{1}+l_{2}n_{2}+\ldots+l_{s}n_{s}=0.
		\end{equation}
		with respect to $s$ variables $l_{1},l_{2},\ldots, l_{s}$.
	Let $\mathcal{W}$ be the solution space of equation (\ref{equ-w-4}). Then $m_{G}(1)=\dim \zeta_{G}(1)=\dim \mathcal{W}$. In fact, $\zeta_{G}(1)= \{\textbf{x}|L(G)\cdot\textbf{x}=1\cdot \textbf{x}\}$ and $\dim\mathcal{W}=s-1=\omega(G-v)-1$.
	Therefore, it deduce that
	\begin{equation*}
		m_{G}(1)=\omega(G-u)-1,
	\end{equation*}
	as desired.
\end{proof}
From Lemma \ref{lem-w-9}, the following corollary is deduced immediately.
\begin{cor}\label{cor-k-1}
	If $G=K_{1}\vee G_{1}$, then $m_{G}(1)=\omega(G_{1})-1$.
\end{cor}

Now we shall give the proof of Theorem \ref{thm-k-5} in the following.
\begin{proof}[\textbf{Proof of Theorem \ref{thm-k-5}}]
	Let $G$ be a $k$-cyclic graph with $m_{G}(1)=t\ge 0$. If $t=0$ then $m_{G}(1)=0\le n-\big\lceil\frac{3+\sqrt{8k+1}} {2}\big\rceil$ since $n\ge n_{k}=\big\lceil\frac{3+\sqrt{8k+1}} {2}\big\rceil$ and $k\geq0$; and otherwise, then by Lemma \ref{lem-2-4}, we have $m_{\overline{G}}(n-1)=m_{G}(1)=t\ge 1$. Furthermore, it follows from Lemma \ref{lem-2-5} that
	\begin{equation}\label{equa-2}
		t(n-1)\leq\mu_{1}(\overline{G})+\mu_{2}(\overline{G})+\ldots+\mu_{t}(\overline{G})\leq \frac{n(n-1)}{2}-(n+k-1)+\frac{t(t+1)}{2}.
	\end{equation}
	Hence, one can obtain that
	\begin{equation}\label{equa-3}
		t^{2}-(2n-3)t+n^{2}-3n+2-2k\geq 0.
	\end{equation}
	The discriminant of the quadratic polynomial in (\ref{equa-3}) is given by
	\begin{align*}
		\Delta
		&=(2n-3)^2-4\bigl(n^2-3n+2-2k\bigr)\\
		&=4n^2-12n+9-4n^2+12n-8+8k\\
		&=8k+1>0 \ \  (\mbox{since} \ \ k\geq0).
	\end{align*}
	Thus, the corresponding quadratic equation has two
	distinct real roots
	\begin{equation*}
		t_{1}=n-\frac{3+\sqrt{8k+1}}{2}
		\quad\text{and}\quad
		t_{2}=n-\frac{3-\sqrt{8k+1}}{2}.
	\end{equation*}
	Since the leading coefficient is positive, it follows from (\ref{equa-3}) that
	\begin{equation*}
		t\leq n-\frac{3+\sqrt{8k+1}}{2}
		\quad\text{or}\quad
		t\geq n-\frac{3-\sqrt{8k+1}}{2}.
	\end{equation*}
	Assume by a contradiction that $t\geq n-\frac{3-\sqrt{8k+1}}{2}$. Since $G$ is a connected graph we have $t \leq n-1$. Thus, $n-\frac{3-\sqrt{8k+1}}{2}\le n-1$. It further deduces that $k \le \frac{1}{4}$. Because $k$ is non-negative integer, this immediately implies $k=0$. Substituting this one can obtain that
	$$\displaystyle t\geq n-\frac{3-\sqrt{8k+1}}{2}=n-1,$$ from which we conclude $t=n-1$. Thus, we have $\operatorname{Spec}_{L}(G)=\{1^{n-1},0^{1}\}$, it implies that  $\sum_{i=1}^{n}\mu_{i}(G)=n-1$. However, $\sum_{i=1}^{n}\mu_{i}(G)=2|E(G)|=2(n-1)$, it therefore follows that $n=1$, this contradicts $n\geq3$ while $k=0$.
	Thus, we have $t\leq n-\frac{3+\sqrt{8k+1}}{2}$, and furthermore, $$m_{G}(1)\leq n-\big\lceil\frac{3+\sqrt{8k+1}} {2}\big\rceil$$ for any integer $k\ge0$. 
	\\
	More precisely, if $\displaystyle k={\widehat{n}_{k}-1\choose 2}$, then $n_{k}=\widehat{n}_{k}$ is an integer, and thus, $$m_{G}(1)\leq n-\frac{3-\sqrt{8k+1}}{2}.$$ If  $\displaystyle {\widehat{n}_{k}-1\choose 2}<k<{\widehat{n}_{k}\choose 2}$, then $n_{k}$ is a fraction. Hence, $$\displaystyle m_{G}(1)\leq  n-\big[\frac{3+\sqrt{8k+1}}{2}\big]-1.$$
	
	Now we turn to verify the equalities, and distinct two cases in the following.
	
	\textbf{Case 1.} $\displaystyle n=n_{k}$.
	
	In this case, $m_{G}(1)=0$ implies that $G$ is a minimal $k$-cyclic graph. Conversely, if $G$ is a minimal $k$-cyclic graph, we assume that $m_{G}(1)\ge 1$. Then from the proof above, $m_{G}(1)\le n-n_{k}$. However, since $G$ is a minimal $k$-cyclic graph, we have $n=n_{k}$ by definition. Substituting this into the above inequality yields $m_{G}(1)\le 0$, which directly contradicts the assumption that $m_{G}(1)\ge 1$. We therefore conclude that $m_{G}(1)=0$.

	\textbf{Case 2.} $\displaystyle n\geq n_{k}+1$.
	
	In this case we also consider two subcases in light of the value of $k$.
	
	\textbf{Case 2.1.} When $\displaystyle k={\widehat{n}_{k}-1\choose 2}$, we suppose that $\mathbb{G}$ is the extreme $k$-cyclic graph with $m_{\mathbb{G}}(1)=n-\frac{3+\sqrt{8k+1}}{2}$. Note that $\widehat{n}_{k}=n_{k}$ in this case. Hence, $m_{\mathbb{G}}(1)=n-\frac{3+\sqrt{8k+1}} {2}=n-n_{k}\ge 1$, which means that all equalities hold in (\ref{equa-2}), and thus, together with Lemma \ref{lem-2-5} we have that $\mu_{1}(\overline{\mathbb{G}})=\mu_{2}(\overline{\mathbb{G}})=\ldots=\mu_{n-n_{k}}(\overline{\mathbb{G}})=n-1$, and $\overline{\mathbb{G}}$ is a threshold graph of order $n$ with clique number $\beta(\overline{\mathbb{G}})=n-n_{k}+1$. Therefore, it follows from Lemma \ref{lem-2-4} and Fact \ref{fact-w-1} that $\mu_{n-1}(\mathbb{G})=\mu_{n-2}(\mathbb{G})=\ldots=\mu_{n_{k}}(\mathbb{G})=1$, and $\mathbb{G}$ is also a threshold graph of order $n$ with $\alpha(\mathbb{G})=n-n_{k}+1$ and $\beta(\mathbb{G})=n_{k}$.
	\begin{claim}\label{claim-w-1}
		$\mathbb{G}$ contains a universal cut vertex $u$ (say).
	\end{claim}
	\begin{proof}
		Assume that $\mathbb{G}$ has no cut vertex, then by Definition \ref{def-w-1} we have $\mathbb{G} \cong K_{n}$. However, the Laplacian spectrum of $K_{n}$ is $\operatorname{Spec}_{L}(K_{n})=\{n^{n-1}, 0^{1}\}$. From $\mu_{n-1}(\mathbb{G})=1=\mu_{n-1}(K_{n})=n$ we have $n=1$, it is a contradiction to $n\ge n_{k}+1\ge3$. Hence, $\mathbb{G}$ has a cut vertex, we notice that $\mu_{n-1}(\mathbb{G}) = 1$, it therefore follows from Lemma \ref{lem-2-6} that the claim holds.
	\end{proof}
	
	\begin{claim}\label{claim-w-2}
		Let $p(\mathbb{G})$ be the number of pendant vertices in $\mathbb{G}$. Then $p(\mathbb{G})\le n-n_{k}+1$. Moreover, if $\mathbb{G}$ is a tree then $p(\mathbb{G})=n-n_{k}+1$, and otherwise $p(\mathbb{G})\le n-n_{k}$.
	\end{claim}
	\begin{proof}
		By Claim \ref{claim-w-1} we see that $\mathbb{G}$ has a universal cut vertex $u$, hence all pendant vertices of $\mathbb{G}$ must be adjacent to $u$, that is, $u$ is the just qusi-pendant vertex of $\mathbb{G}$, and this moment $q(\mathbb{G})=1$. Furthermore, by Lemma \ref{lem-2-7} we have $n-n_{k} =m_{\mathbb{G}}(1)\ge p(\mathbb{G})-1$, and so, $p(\mathbb{G})\le n-n_{k}+1$.
		
		Especially, if $\mathbb{G}$ is a tree (i.e., $k=0$), then by Lemma \ref{lem-2-8} we have $n-n_{k}=m_{\mathbb{G}}(1)\le p(\mathbb{G})-1$, that is, $p(\mathbb{G})\ge n-n_{k}+1$. On the other hand, $p(\mathbb{G})\le n-n_{k}+1$. Therefore, $p(\mathbb{G})= n-n_{k}+1$ (In fact, $n_{k}=2$ due to $k=0$. At this point, $p(\mathbb{G})= n-1$. Noting that $\mathbb{G}$ is an extremal tree with a universal cut vertex $u$, and thus $\mathbb{G}\cong S_{n}$). Otherwise, assume that $p(\mathbb{G})=n-n_{k}+1$, then $\mathbb{G}$ has an independent set of order $n-n_{k}+2$ since $\mathbb{G}$ has a clique of order $n_{k}-1\ge 2$ containing the universal cut vertex $u$. Hence, it contradicts $\alpha(\mathbb{G})=n-n_{k}+1$, and so, $p(\mathbb{G})\le n-n_{k}$. Furthermore, we notice that $\alpha(\mathbb{G})=n-n_{k}+1$ and $\beta(\mathbb{G})=n_{k}$, it thus follows that $P(\mathbb{G})=n-n_{k}$, as required.
	\end{proof}
	Combining with Claims \ref{claim-w-1} and \ref{claim-w-2}, we conclude that $\mathbb{G}\cong K_{1} \vee (K_{n_{k}-1} \cup (n-n_{k})K_{1})$.

	\textbf{Case 2.2.} When $\displaystyle {\widehat{n}_{k}-1\choose 2}<k<{\widehat{n}_{k}\choose 2}$, we assume that $\mathbb{H}$ is the extremal graph with $m_{\mathbb{H}}(1)=n-\big\lceil\frac{3+\sqrt{8k+1}} {2}\big\rceil=n-n_{k}=t\ge 1$, and $\overline{\mathbb{H}}$ be the complement graph. Then by Lemma \ref{lem-2-4}, we have $m_{\overline{\mathbb{H}}}(n-1)=t$ due to $m_{\mathbb{H}}(1)=t$. From the previous discussion of Case 2.1, if $k={\widehat{n}_{k}-1\choose 2}$ then
		\begin{equation}\label{equa-4}
			t(n-1)=\mu_{1}(\overline{\mathbb{G}})+\mu_{2}(\overline{\mathbb{G}})+\ldots+\mu_{t}(\overline{\mathbb{G}})= \frac{n(n-1)}{2}-(n+k-1)+\frac{t(t+1)}{2}.
		\end{equation}
		But in this case,
		\begin{equation}\label{equa-5}
			t(n-1)=\mu_{t}(\overline{\mathbb{H}})+\mu_{t-1}(\overline{\mathbb{H}})+\ldots+\mu_{1}(\overline{\mathbb{H}})< \frac{n(n-1)}{2}-(n+k-1)+\frac{t(t+1)}{2}.
		\end{equation}
		Through a comparison with (\ref{equa-4}) one can observe that $n$ and $t$ are fixed in (\ref{equa-5}), if some edge(s) are deleted from $\overline{\mathbb{H}}$ then the inequality on the right-hand side of (\ref{equa-5}) becomes an equality. In this case, let $\varepsilon(k)$ denote the number of edges deleted from $\overline{\mathbb{H}}$. Then we have
		\begin{equation}\label{equa-6}
			\begin{array}{ll}
				t(n-1)=\mu_{t}(\overline{\mathbb{H}})+\mu_{t-1}(\overline{\mathbb{H}})+\ldots+\mu_{1}(\overline{\mathbb{H}})
				= \frac{n(n-1)}{2}-(n+k-1)-\varepsilon(k)+\frac{t(t+1)}{2}.
			\end{array}
		\end{equation}
		At this moment, $t=n-n_{k}$ and $n_{k}=\big\lceil\frac{3+\sqrt{8k+1}} {2}\big\rceil$. Substituting these into (\ref{equa-6}) and solving for $\varepsilon(k)$, we can obtain
		\begin{equation}\label{equa-7}
			\begin{array}{ll}
				\varepsilon(k) &=\displaystyle \frac{1}{2} \left( \left\lceil \frac{3 + \sqrt{8k + 1}}{2} \right\rceil - 1 \right) \left( \left\lceil \frac{3 + \sqrt{8k + 1}}{2} \right\rceil - 2 \right) - k\\
				&\displaystyle ={n_{k}-1\choose 2}-k.
			\end{array}
		\end{equation}
	Furthermore, by the proof of Case 2.1. we see that, all the equalities of (\ref{equa-6}) hold if and only if
		\begin{equation*}
			\overline{\mathbb{H}}-\varepsilon(k)\cdot e \cong \overline{K_{1} \vee (K_{n_{k}-1} \cup (n-n_{k})K_{1})},
		\end{equation*}
		and so,
		\begin{equation}\label{equa-w-7}
			\mathbb{H} \cong K_{1} \vee (K_{n_{k}-1} \cup (n-n_{k})K_{1})-\varepsilon(k)\cdot e
		\end{equation}
		where $\varepsilon(k)={n_{k}-1\choose 2}-k$.
	
	In fact, from the equalities of the left-side of (\ref{equa-6}) one can see that, $\mathbb{H}$ must contain a universal vertex. Thus, we have
		\begin{equation*}
			K_{1} \vee \big(K_{n_{k}-1} \cup (n-n_{k})K_{1}\big)-\varepsilon(k)\cdot e \cong K_{1} \vee \bigg(\big(K_{n_{k}-1}-\varepsilon(k)\cdot e\big) \cup (n-n_{k})K_{1}\bigg).
		\end{equation*}
		Let $H=K_{n_{k}-1}-\varepsilon(k)\cdot e$. Then by (\ref{equa-w-7}) we get
		\begin{equation*}
			\mathbb{H} \cong K_{1} \vee \bigg(H \cup (n-n_{k})K_{1}\bigg),
		\end{equation*}
		here the induced subgraph $K_{1} \vee H$ of $\mathbb{H}$ is a minimal $k$-cyclic graph, and $H$ is of order $n_{k}-1$ and size $k$.
	
	In order to characterize an exception, we give a claim as follows.
		\begin{claim}\label{claim-w-3}
			$1\leq \varepsilon(k)\leq n_{k}-3$
		\end{claim}
		\begin{proof}
			Since $\displaystyle {\widehat{n}_{k}-1\choose 2}<k<{\widehat{n}_{k}\choose 2}$, and $n_{k}=\widehat{n}_{k}+1$, we substitute these relations into (\ref{equa-7}) and obtain
			\begin{equation*}
				\displaystyle {{n}_{k}-1\choose 2}-{\widehat{n}_{k}\choose 2}< \varepsilon(k)< {{n}_{k}-1\choose 2}-{\widehat{n}_{k}-1\choose 2},
			\end{equation*}
			and hence,
			\begin{equation*}
				\displaystyle {\widehat{n}_{k}\choose 2}-{\widehat{n}_{k}\choose 2}< \varepsilon(k)< {\widehat{n}_{k}\choose 2}-{\widehat{n}_{k}-1\choose 2},
			\end{equation*}
			it therefore deduces that
			\begin{equation*}
				\displaystyle 0<  \varepsilon(k)< \widehat{n}_{k}-1.
			\end{equation*}
			i.e., $1\leq \varepsilon(k)\leq \widehat{n}_{k}-2$. Thus, $1\leq \varepsilon(k)\leq n_{k}-3$.
	\end{proof}
	
	For any $\varepsilon(k)$, let $\mathbb{H}_{1}=K_{1} \vee\bigg(\big(K_{n_{k}-2} \cup (n-n_{k}+1)K_{1}\big)+(n_{k}-2-\varepsilon(k))\cdot e\bigg)$ by Claim \ref{claim-w-3}.
		Since $\mathbb{H}=K_{1} \vee \bigg(\big(K_{n_{k}-1}-\varepsilon(k)\cdot e\big) \cup (n-n_{k})K_{1}\bigg)$ is a $k$-cyclic graph and $\displaystyle {\widehat{n}_{k}-1\choose 2}<k<{\widehat{n}_{k}\choose 2}$,
		both $\mathbb{H}$ and $\mathbb{H}_{1}$ have the same cyclomatic number.
		From Corollary \ref{cor-k-1} one can see that, if  $m_{\mathbb{H}}(1)=m_{\mathbb{H}_{1}}(1)$ then $n_{k}-2-\varepsilon(k)$ edges are only added between some pendant vertices of $\mathbb{H}_{1}$, say $s$ such pendant vertices. Then we have
		\begin{equation*}
			\mathbb{H}_{1}=K_{1} \vee\big(K_{n_{k}-2} \cup (n-n_{k}+1-s)K_{1}\cup H'_{1}\big)
		\end{equation*}
		where $H'_{1}$ is an induced connected subgraph of $s$ pendant vertices as well as $n_{k}-2-\varepsilon(k)$ edges in $\mathbb{H}_{1}$. In this case, $m_{\mathbb{H}_{1}}(1)= m_{\mathbb{H}}(1)$ together with Corollary \ref{cor-k-1}   we have
		\begin{equation*}
			n-n_{k}+1-s+2-1=n-n_{k},
		\end{equation*}
		and so, it follows that $s=2$, which implies that $n_{k}-2-\varepsilon(k)=1$, that is, $\varepsilon(k)=n_{k}-3$. Hence,
		\begin{equation*}
			\mathbb{H}_{1}=K_{1} \vee\big(K_{n_{k}-2} \cup K_{2}\cup (n-n_{k}-1)K_{1}\big)
	\end{equation*}
	
	Consequently, if $\varepsilon(k)=n_{k}-3$, both two graphs $K_{1} \vee \bigg(\big(K_{n_{k}-1}-(n_{k}-3)\cdot e\big) \cup (n-n_{k})K_{1}\bigg)$ and $K_{1} \vee\big(K_{n_{k}-2} \cup K_{2}\cup (n-n_{k}-1)K_{1}\big)$ has the same multiplicity of the Lapalcian eigenvalue $1$. 
	
	The proof is complete.
\end{proof}

%\begin{cor}\label{thm-w-1}
%For a fixed $k\ge 0$, let $G$ be a minimal $k$-cyclic graph. Then $m_{G}(1)=0$.
%\end{cor}

\section{Concluding remarks}\label{sec-4}
In this paper, we first prove that every $k$-cyclic graph with fixed $k\geq 0$ is of order $n\geq \big\lceil\frac{3+\sqrt{8k+1}} {2}\big\rceil$, see Proposition \ref{prop-k-1}. After then, we apply the Brouwer's inequality to obtain a sharp upper bound $m_{G}(1)\leq  n-\big\lceil\frac{3+\sqrt{8k+1}} {2}\big\rceil$ for $m_{G}(1)$ of a $k$-cyclic graph $G$ with fixed order $n$ and $k\geq0$, see Theorem \ref{thm-k-5}. Additionally, we completely determine all extremal graphs for which equality holds. Furthermore, our results not only generalize the main results of Guo et al.\cite{Guo}, Wen et al.\cite{Wen}, and Wang et al.\cite{Wang}, but also answer the question posed in Problem \ref{prob-1} concerning which graphs achieve the maximum multiplicity of the Laplacian eigenvalue $1$.

It is well known that every connected graph can be viewed as a $k$-cyclic graph with $k=\varepsilon-n+1$. Therefore, our results also apply to any connected graph. In addition, the multiplicity of the Laplacian eigenvalue $1$ of a graph equals the sum of that of its components. Consequently, our results remain valid for any simple graphs.

Given a graph family $\mathcal{H}$. A graph is said to be \emph{$\mathcal{H}$-free} if it does not contain subgraphs isomorphic to any members of $H$. A natural question is that if $G$ is $\mathcal{H}$-free, what is the maximum multiplicity of the Laplacian eigenvalue $1$ of $G$ ?  Hence, we end this paper by the following problem for future research.
\begin{prob}\label{prob-2}
	Let $G$ be a $\mathcal{H}$-free graph with fixed order $n\ge 2$. Determine the maximum multiplicity of the Laplacian eigenvalue $1$ of the graph $G$, and then characterized the corresponding extremal graphs.
\end{prob}

\section*{Data availability statement}
Our manuscript has no associated data.
\section*{Declaration of competing interest}
The authors declare that they have no conflicts of interest to this work.

\section*{Acknowledgements}
This paper relies heavily on the main conclusions established in two papers by   Kothari et al.\cite{Kothari} and  Cai et al.\cite{Cai} independently, in which Kothari et al.\cite{Kothari} proved the full Brouwer's Laplacian spectrum conjecture, and Cai et al.\cite{Cai} completely characterized all graphs with its equality. The present author would like to thank them for their papers. Without those, this one would not have been written.

The first author (F. Wen) was partially supported by the National Natural Science Foundation of China (Nos. 11961041; 12261055), the Key Project of Natural Science Foundation of Gansu Province, China (No. 24JRRA222) and the Foundation for Innovative Fundamental Research Group Project of Gansu Province, China (No. 25JRRA805); the third author (T. Wang) was partially supported by the National Natural Science Foundation of China (No. 12601667).

\end{document}